\documentclass[10pt]{amsart}
\usepackage{amsmath,amssymb,amstext,amsthm,amscd}

\theoremstyle{plain}

\newtheorem*{lemma*}{Lemma}

\newtheorem*{theorem*}{Theorem}

\newtheorem*{claim*}{Claim}

\newtheorem*{conjecture*}{Conjecture}
\newtheorem*{question*}{Question}

\theoremstyle{definition}

\newtheorem*{definition*}{Definition}

\theoremstyle{remark}
\newtheorem*{remark*}{Remark}
\newtheorem*{remarks*}{Remarks}
\newtheorem*{example*}{Example}
\newtheorem*{exercise*}{Exercise}
\newtheorem*{acknowledgement*}{Acknowledgement}

\newcommand{\pp}{\mathbb{P}}

\newcommand{\Mbar}{\overline{\mathcal{M}}}

\def\({\left(}
\def\){\right)}
\def\<{\langle}
\def\>{\rangle}
\def\d{\partial}

\def\cM{{\mathcal M}}
\def\ocM{{\overline{\cM}}}

\thanks{The second author is partially supported by NSF.}

\begin{document}

\title{Tautological equations in $\ocM_{3,1}$ via invariance conjectures}

\author{D.~Arcara}
\address{Department of Mathematics, University of Utah,
155 S. 1400 E., Room 233, Salt Lake City, UT 84112-0090, USA}
\email{arcara@math.utah.edu}

\author{Y.-P.~Lee}
\email{yplee@math.utah.edu}

\maketitle

\begin{center}
\emph{This paper is dedicated to Matteo S.~Arcara, \\
who was born while the paper was at its last phase of preparation.}
\end{center}

\begin{abstract}
A new tautological equation of $\ocM_{3,1}$ in codimension 3
is derived and proved, using the invariance condition explained in
\cite{AL, ypL1, ypL2, ypL3}.
\end{abstract}

\section{Introduction}

This work is a continuation of \cite{AL}. We apply the same technique 
to $R^3(\Mbar_{3,1})$ to find a tautological equation.
A general scheme and practical steps, as well as notations used
in this paper, can be found in \cite{ypL3} and \cite{AL}.

Our way of finding this equation is fairly simple.
\begin{enumerate}
\item By the Hodge number calculations done by E.~Getzler and  E.~Looijenga
\cite{GL}, there is any new equation in $R^3(\Mbar_{3,1})$.
\item Apply \emph{Invariance Conjecture~1} (Theorem~5 \cite{ypL3}) to obtain
the coefficients of the equation.
\end{enumerate}

(2) gives a \emph{necessary} condition.
Combined with (1), this generates and proves the new equation.

Our motivation was quite simple. In the earlier work \cite{GL, AL},
the conjectural framework is shown to be valid in genus one and two.
While it is satisfactory to learn that all previous equations can be 
derived in this framework, the Conjectures predicts the possibility of 
finding all tautological equations under this framework.
This work is set to be the first step towards this goal.

The choice of codimension 3 in $\Mbar_{3,1}$ is almost obvious.
First of all, the Invariance Conjectures works inductively.
Given what we know about genus one and two, it is only reasonable to proceed
to either $\Mbar_{2,n}$ for $n \ge 4$ or $\Mbar_{3,1}$. 
Secondly, one also knows from the Theorem of Graber--Vakil \cite{GV} that 
$\psi^3$ on $\Mbar_{3,1}$ is rational equivalent to a sum of boundary strata
containing at least one (geometrical) genus zero component. 
Thirdly, Getzler and Looijenga \cite{GL} have shown that there is 
only one relation in codimension 3 in $\Mbar_{3,1}$.
That makes it a reasonable place to start. 

The main result of this paper is the following:
\begin{theorem*}
There is a new tautological equation for codimension 3 strata in 
$\Mbar_{3,1}$.

\begin{equation*}
 \begin{split}
\< \d^x_3 \>_3 =
 &\frac{5}{72} \< \d^x \d^{\mu_1} \d^{\mu_2} \> \< \d^{\mu_1} \d^{\mu_3} \d^{\mu_2} \> \< \d^{\mu_3} \>_2
+ \frac{1}{252} \< \d^{\mu_1} \d^{\mu_2} \d^{\mu_2} \> \< \d^x_1 \d^{\mu_1} \>_2 \\
+& \frac{5}{72} \< \d^x \d^{\mu_1} \d^{\mu_1} \d^{\mu_2} \> \< \d^{\mu_2}_1 \>_2 
+ \frac{5}{42} \< \d^x \d^{\mu_1} \d^{\mu_2} \> \< \d^{\mu_1}_1 \d^{\mu_2} \>_2 
+ \frac{41}{21} \< \d^x \d^{\mu_1} \d^{\mu_2} \> \< \d^{\mu_1} \>_1 \< \d^{\mu_2}_1 \>_2 \\
+& \frac{11}{40320} \< \d^x \d^{\mu_1} \d^{\mu_2} \d^{\mu_2} \> \< \d^{\mu_1} \d^{\mu_3} \d^{\mu_3} \>_1 
+ \frac{1}{13440} \< \d^x \d^{\mu_1} \d^{\mu_2} \d^{\mu_2} \>_1 \< \d^{\mu_1} \d^{\mu_3} \d^{\mu_3} \> \\
+& \frac{1}{8064} \< \d^x \d^{\mu_1} \>_1 \< \d^{\mu_1} \d^{\mu_2} \d^{\mu_2} \d^{\mu_3} \d^{\mu_3} \> 
+ \frac{191}{120960} \< \d^x \d^{\mu_1} \d^{\mu_2} \d^{\mu_2} \d^{\mu_3} \d^{\mu_3} \> \< \d^{\mu_1} \>_1 \\
+& \frac{1}{5040} \< \d^x \d^{\mu_1} \d^{\mu_2} \> \< \d^{\mu_1} \d^{\mu_2} \d^{\mu_3} \d^{\mu_3} \>_1 
+ \frac{1}{4032} \< \d^x \d^{\mu_1} \d^{\mu_2} \>_1 \< \d^{\mu_1} \d^{\mu_2} \d^{\mu_3} \d^{\mu_3} \> \\
+& \frac{17}{2880} \< \d^x \d^{\mu_1} \d^{\mu_2} \d^{\mu_3} \d^{\mu_3} \> \< \d^{\mu_1} \d^{\mu_2} \>_1 
+ \frac{1}{840} \< \d^x \d^{\mu_1} \>_1 \< \d^{\mu_1} \d^{\mu_2} \>_1 \< \d^{\mu_2} \d^{\mu_3} \d^{\mu_3} \> \\
+& \frac{1}{336} \< \d^x \d^{\mu_1} \d^{\mu_2} \>_1 \< \d^{\mu_1} \>_1 \< \d^{\mu_2} \d^{\mu_3} \d^{\mu_3} \> 
+ \frac{1}{126} \< \d^x \d^{\mu_1} \d^{\mu_2} \> \< \d^{\mu_1} \>_1 \< \d^{\mu_2} \d^{\mu_3} \d^{\mu_3} \>_1 \\
+& \frac{23}{5040} \< \d^x \d^{\mu_1} \d^{\mu_2} \d^{\mu_2} \> \< \d^{\mu_1} \d^{\mu_3} \>_1 \< \d^{\mu_3} \>_1
+ \frac{17}{5040} \< \d^x \d^{\mu_1} \>_1 \< \d^{\mu_1} \d^{\mu_2} \d^{\mu_3} \d^{\mu_3} \> \< \d^{\mu_2} \>_1 \\
+& \frac{113}{2520} \< \d^x \d^{\mu_1} \d^{\mu_2} \d^{\mu_3} \d^{\mu_3} \> \< \d^{\mu_1} \>_1 \< \d^{\mu_2} \>_1+
 \frac{1}{210} \< \d^x \d^{\mu_1} \d^{\mu_2} \d^{\mu_3} \> \< \d^{\mu_1} \d^{\mu_2} \d^{\mu_3} \>_1\\
+& 0 \< \d^x \d^{\mu_1} \d^{\mu_2} \d^{\mu_3} \>_1 \< \d^{\mu_1} \d^{\mu_2} \d^{\mu_3} \> +
 \frac{1}{84} \< \d^x \d^{\mu_1} \>_1 \< \d^{\mu_1} \d^{\mu_2} \d^{\mu_3} \> \< \d^{\mu_2} \d^{\mu_3} \>_1\\
+& \frac{211}{1260} \< \d^x \d^{\mu_1} \d^{\mu_2} \d^{\mu_3} \> \< \d^{\mu_1} \d^{\mu_2} \>_1 \< \d^{\mu_3} \>_1+
 \frac{1}{1260} \< \d^x \d^{\mu_1} \d^{\mu_2} \>_1 \< \d^{\mu_1} \d^{\mu_2} \d^{\mu_3} \> \< \d^{\mu_3} \>_1 \\
+& \frac{1}{630} \< \d^x \d^{\mu_1} \d^{\mu_2} \> \< \d^{\mu_1} \d^{\mu_2} \d^{\mu_3} \>_1 \< \d^{\mu_3} \>_1+
 \frac{11}{140} \< \d^x \d^{\mu_1} \d^{\mu_2} \> \< \d^{\mu_1} \d^{\mu_3} \>_1 \< \d^{\mu_2} \d^{\mu_3} \>_1\\
+& \frac{4}{35} \< \d^x \d^{\mu_1} \d^{\mu_2} \> \< \d^{\mu_1} \>_1 \< \d^{\mu_2} \d^{\mu_3} \>_1 \< \d^{\mu_3} \>_1 +
 \frac{2}{105} \< \d^x \d^{\mu_1} \>_1 \< \d^{\mu_1} \d^{\mu_2} \d^{\mu_3} \> \< \d^{\mu_2} \>_1 \< \d^{\mu_3} \>_1\\
+& \frac{89}{210} \< \d^x \d^{\mu_1} \d^{\mu_2} \d^{\mu_3} \> \< \d^{\mu_1} \>_1 \< \d^{\mu_2} \>_1 \< \d^{\mu_3} \>_1+
 \frac{1}{53760} \< \d^x \d^{\mu_1} \d^{\mu_1} \d^{\mu_2} \d^{\mu_2} \d^{\mu_3} \d^{\mu_3} \>.
 \end{split}
\end{equation*}
\end{theorem*}

Note that the notations employed here are called gwis, defined in \cite{ypL3}.
They are a convenient way of presenting decorated graphs.

\begin{remarks*}
(i) While our paper was under preparation, a preprint by T.~Kimura and X.~Liu 
\cite{KL} appeared on the arxiv. There are two major differences between
our results. First, their choice of basis of codimension 3 strata in
$\Mbar_{3,1}$ is different. They use
$(3') := \< \d^{\mu_1} \d^{\mu_2} \d^{\mu_2} \> \< \d^x \d^{\mu_1}_1 \>_2 $
instead of (3) below. 
We have checked that their equation is equivalent to ours.
Second, their approach is ``traditional'': knowing there must be a
relation from Graber--Vakil, they can then proceed to find the coefficients
based on the evaluation of the Gromov--Witten invariants of $\pp^1$.

Our approach is quite different.
There are no computer-aided calculation of the Gromov--Witten invariants.
Only linear algebra is involved in the calculation.

(ii) Using the same technique, the reader may amuse himself with the following
exercises.

\emph{There is no (new) relation among the boundary strata and $\psi$ classes 
in codimension 2 in $\Mbar_{3,1}$.}

This can be shown, for example, by the same technique used in this paper.
Let $(i)$ denotes a basis of these strata.
When one sets a hypothetical equation $c_i (i) =0$ and 
imposes the invariance condition, the only solution is $c_i=0$ for all $i$.
\end{remarks*}

\begin{acknowledgement*}
We are thankful to E.~Getzler and R.~Vakil for the useful discussions.
Part of the work was done during the second author's stay in NCTS, whose
hospitality is greatly appreciated.
\end{acknowledgement*}

\section{Strata of $\ocM_{3,1}$}

We start with enumerating codimension 3 strata in $\Mbar_{3,1}$.

Out of the several strata (allowing $\psi$ classes) of codimension $3$ in
$\ocM_{3,1}$, many of them can be written in terms of the others using WDVV,
TRR's, Mumford--Getzler's, Getzler's and Belorousski-Pandharipande equations.
After applying those equations, we can write all of the strata in terms of the
following ones:

$$ \begin{matrix}
& (1) & \< \d^x_3 \>_3
& (2) & \< \d^x \d^\alpha \d^\beta \> \< \d^\mu \d^\alpha \d^\beta \> \< \d^\mu \>_2 \\
& (3) & \< \d^\mu \d^\alpha \d^\alpha \> \< \d^x_1 \d^\mu \>_2 
& (4) & \< \d^x \d^\mu \d^\mu \d^\nu \> \< \d^\nu_1 \>_2 \\
& (5) & \< \d^x \d^\mu \d^\nu \> \< \d^\mu_1 \d^\nu \>_2
& (6) & \< \d^x \d^\mu \d^\nu \> \< \d^\mu \>_1 \< \d^\nu_1 \>_2 \\
& (7) & \< \d^x \d^\mu \d^\nu \d^\nu \> \< \d^\mu \d^\alpha \d^\alpha \>_1
& (8) & \< \d^x \d^\mu \d^\nu \d^\nu \>_1 \< \d^\mu \d^\alpha \d^\alpha \> \\
& (9) & \< \d^x \d^\mu \>_1 \< \d^\mu \d^\nu \d^\nu \d^\alpha \d^\alpha \> 
& (10) & \< \d^x \d^\mu \d^\nu \d^\nu \d^\alpha \d^\alpha \> \< \d^\mu \>_1 \\
& (11) & \< \d^x \d^\mu \d^\nu \> \< \d^\mu \d^\nu \d^\alpha \d^\alpha \>_1
& (12) & \< \d^x \d^\mu \d^\nu \>_1 \< \d^\mu \d^\nu \d^\alpha \d^\alpha \> \\
& (13) & \< \d^x \d^\mu \d^\nu \d^\alpha \d^\alpha \> \< \d^\mu \d^\nu \>_1
& (14) & \< \d^x \d^\mu \>_1 \< \d^\mu \d^\nu \>_1 \< \d^\nu \d^\alpha \d^\alpha \> \\
& (15) & \< \d^x \d^\mu \d^\nu \>_1 \< \d^\mu \>_1 \< \d^\nu \d^\alpha \d^\alpha \> 
& (16) & \< \d^x \d^\mu \d^\nu \> \< \d^\mu \>_1 \< \d^\nu \d^\alpha \d^\alpha \>_1 \\
& (17) & \< \d^x \d^\mu \d^\nu \d^\nu \> \< \d^\mu \d^\alpha \>_1 \< \d^\alpha \>_1
& (18) & \< \d^x \d^\mu \>_1 \< \d^\mu \d^\nu \d^\alpha \d^\alpha \> \< \d^\nu \>_1 \\
& (19) & \< \d^x \d^\mu \d^\nu \d^\alpha \d^\alpha \> \< \d^\mu \>_1 \< \d^\nu \>_1
& (20) & \< \d^x \d^\mu \d^\nu \d^\alpha \> \< \d^\mu \d^\nu \d^\alpha \>_1 \\
& (21) & \< \d^x \d^\mu \d^\nu \d^\alpha \>_1 \< \d^\mu \d^\nu \d^\alpha \> 
& (22) & \< \d^x \d^\mu \>_1 \< \d^\mu \d^\nu \d^\alpha \> \< \d^\nu \d^\alpha \>_1 \\
& (23) & \< \d^x \d^\mu \d^\nu \d^\alpha \> \< \d^\mu \d^\nu \>_1 \< \d^\alpha \>_1
& (24) & \< \d^x \d^\mu \d^\nu \>_1 \< \d^\mu \d^\nu \d^\alpha \> \< \d^\alpha \>_1 \\
& (25) & \< \d^x \d^\mu \d^\nu \> \< \d^\mu \d^\nu \d^\alpha \>_1 \< \d^\alpha \>_1
& (26) & \< \d^x \d^\mu \d^\nu \> \< \d^\mu \d^\alpha \>_1 \< \d^\nu \d^\alpha \>_1 \\
& (27) & \< \d^x \d^\mu \d^\nu \> \< \d^\mu \>_1 \< \d^\nu \d^\alpha \>_1 \< \d^\alpha \>_1 
& (28) & \< \d^x \d^\mu \>_1 \< \d^\mu \d^\nu \d^\alpha \> \< \d^\nu \>_1 \< \d^\alpha \>_1 \\
& (29) & \< \d^x \d^\mu \d^\nu \d^\alpha \> \< \d^\mu \>_1 \< \d^\nu \>_1 \< \d^\alpha \>_1
& (30) & \< \d^x \d^\mu \d^\mu \d^\nu \d^\nu \d^\alpha \d^\alpha \>
\end{matrix} $$

\section{Setting $\mathfrak{r}_1(E)=0$}

Let $E$ be a generic linear combination of these strata
\[
 E := \sum_{k=1}^{30} c_k (k),
\]
where $c_k$ are variables with values in $\mathbb{Q}$.
The Invariance Conjectures predict
\[
  \mathfrak{r}_1 (E)=0.
\]
For the output graphs of $\mathfrak{r}_1(E)$,
we will pick a basis for the tautological algebra and set 
the coefficients of the basis equal to 0. 
This will produce a system of linear equations on $c_k$, (1)-(49), which
then determines $c_k$ completely. 
The final equation $E=0$, with the specified coefficients, is thus obtained.

Note that the (new) half-edges $i,j$ in the output graphs 
are always assumed to be symmetrized.
Some of the graphs will be disconnected.
They are easier to deal with as they involve less relations (e.g.~WDVV).
So let us start with the disconnected terms. 
In each equation, the first column is a basis vector, followed by
its coefficient which is set to zero.

\begin{equation} \label{eq1}
 \< \d^x \d^\alpha \d^\beta \> \< \d^j \d^\alpha \d^\beta \> \< \d^i_1 \>_2
:\quad 
c_2 + c_4 = 0.
\end{equation}
\begin{equation}\label{eq2}
\< \d^x \d^i \d^\alpha \> \< \d^j \d^\beta \d^\beta \> \< \d^\alpha_1 \>_2
: \quad
3 c_3 - c_4 + \frac{1}{24} c_6 = 0.
\end{equation}
\begin{equation}\label{eq3}
\< \d^j \d^\nu \d^\nu \> \< \d^\alpha \d^\alpha \d^\beta \> \< \d^x \d^i \d^\beta \>_1
: \quad
- \frac{1}{80} c_3 - c_8 + \frac{1}{24} c_{15} = 0.
\end{equation}
\begin{equation}\label{eq4}
\< \d^i \d^\nu \d^\beta \> \< \d^\beta \d^x \d^\nu \> \< \d^j \d^\alpha \d^\alpha \>_1
: \quad
c_7 - c_8 - c_{11} = 0.
\end{equation}
\begin{equation}\label{eq5}
 \< \d^j \d^\nu \d^\nu \> \< \d^x \d^\alpha \d^\beta \> \< \d^i \d^\alpha \d^\beta \>_1 
: \quad
\frac{1}{30} c_3 - c_{11} + \frac{1}{24} c_{25} = 0.
\end{equation}
\begin{equation}\label{eq6}
 \< \d^j \d^\nu \d^\nu \> \< \d^i \d^\alpha \d^\beta \> \< \d^x \d^\alpha \d^\beta \>_1 
: \quad
\frac{1}{30} c_3 + 2 c_8 - c_{12} +\frac{1}{24}c_{24}=0.
\end{equation}
\begin{equation}\label{eq7}
 \< \d^j \d^\nu \d^\nu \> \< \d^x \d^i \d^\alpha \> \< \d^\alpha \d^\beta \d^\beta \>_1 
: \quad
- \frac{1}{30} c_3 - c_7 + c_8 + \frac{1}{24} c_{16} = 0.
\end{equation}
\begin{equation}\label{eqnew2}
 \< \d^\mu \d^\nu \d^\nu \> \< \d^j \d^\mu \d^\alpha \> \< \d^i \d^x \>_1 \< \d^\alpha \>_1
: \quad
c_{14} - c_{15} + c_{18} - c_{24} = 0;
\end{equation}
\begin{equation}\label{eq8}
 \< \d^\nu \d^\alpha \d^\alpha \> \< \d^i \d^x \d^\nu \> \< \d^j \d^\mu \>_1 \< \d^\mu \>_1 
: \quad
c_{15} - c_{17} + c_{25} = 0.
\end{equation}
\begin{equation}\label{eq9}
 \< \d^j \d^\nu \d^\nu \> \< \d^x \d^\alpha \d^\beta \> \< \d^i \d^\alpha \>_1 \< \d^\beta \>_1
: \quad
\frac{4}{5} c_3 - c_{16} + \frac{1}{24} c_{27} = 0.
\end{equation}
\begin{equation}\label{eqnew1}
 \< \d^\mu \d^\nu \d^\nu \> \< \d^x \d^i \d^\alpha \> \< \d^j \d^\mu \>_1 \< \d^\alpha \>_1 
: \quad
- c_3 - \frac{1}{240} c_6 + c_{14} - c_{15} = 0.
\end{equation}
\begin{equation}\label{eq10}
 \< \d^j \d^\nu \d^\nu \> \< \d^i \d^\alpha \d^\beta \> \< \d^x \d^\alpha \>_1 \< \d^\beta \>_1
: \quad
\frac{4}{5} c_3 + c_{14} + c_{15} - c_{18}
+ \frac{1}{12} c_{28} = 0.
\end{equation}
\begin{equation}\label{eq11}
\< \d^j \d^\nu \d^\nu \> \< \d^x \d^i \d^\alpha \> \< \d^\alpha \d^\beta \>_1 \< \d^\beta \>_1 
: \quad
- \frac{4}{5} c_3 + c_{15} - c_{17} + \frac{1}{24} c_{27} = 0.
\end{equation}
\begin{equation}\label{eq12}
 \< \d^x \d^j \d^\mu \> \< \d^\mu \>_1 \< \d^i \d^\alpha \d^\beta \> \< \d^\alpha \d^\beta \>_1 
: \quad
\frac{1}{10} c_6 + 2 c_{16} + c_{22} - c_{23} = 0.
\end{equation}
\begin{equation}\label{eq13}
\< \d^x \d^\mu \d^j \> \< \d^i \d^\alpha \d^\beta \> \< \d^\alpha \>_1 \< \d^\beta \>_1 \< \d^\mu \>_1 
: \quad
\frac{7}{10} c_6 + c_{27} + c_{28} - 3 c_{29} = 0.
\end{equation}
\begin{equation}\label{eq14}
 \< \d^j \d^\mu \d^\nu \d^\alpha \> \< \d^\mu \d^\nu \d^\alpha \> \< \d^i \d^x \>_1 = 
\< \d^j \d^\mu \d^\alpha \d^\alpha \> \< \d^\mu \d^\nu \d^\nu \> \< \d^i \d^x \>_1
: \quad
- c_8 + c_9 + \frac{1}{24} c_{14} - c_{21} = 0.
\end{equation}
\begin{equation}\label{eq15}
 \< \d^j \d^\mu \d^\alpha \> \< \d^\mu \d^\nu \d^\nu \d^\alpha \> \< \d^i \d^x \>_1 
: \quad
2 c_9 - c_{12} = 0.
\end{equation}
\begin{equation}\label{eq16}
 \< \d^j \d^\nu \d^\nu \> \< \d^x \d^\alpha \d^\alpha \d^\beta \> \< \d^i \d^\beta \>_1 
: \quad
\frac{1}{48} c_3 - c_7 + \frac{1}{24} c_{17} = 0.
\end{equation}
\begin{equation}\label{eq17}
 \< \d^x \d^i \d^\alpha \d^\alpha \> \< \d^\mu \d^\nu \d^\nu \> \< \d^j \d^\mu \>_1 
: \quad
- \frac{1}{24} c_3 - \frac{1}{240} c_4 - c_8 + \frac{1}{24} c_{14} = 0.
\end{equation}
\begin{equation}\label{eq18}
 \< \d^j \d^\nu \d^\nu \> \< \d^i \d^\alpha \d^\alpha \d^\beta \> \< \d^x \d^\beta \>_1
: \quad
\frac{1}{48} c_3 + c_8 - 2 c_9 + \frac{1}{24} c_{14}
+ \frac{1}{24} c_{18} = 0.
\end{equation}
\begin{equation}\label{eq19}
 \< \d^j \d^\nu \d^\nu \> \< \d^x \d^i \d^\alpha \d^\beta \> \< \d^\alpha \d^\beta \>_1 
: \quad
\frac{7}{30} c_3 + 2 c_8 - c_{13} +\frac{1}{24}c_{23}=0.
\end{equation}
\begin{equation}\label{eq20}
 \< \d^j \d^x \d^\mu \d^\mu \> \< \d^i \d^\alpha \d^\beta \> \< \d^\alpha \d^\beta \>_1 
: \quad
\frac{1}{10} c_4 + 2 c_7 - c_{13} + \frac{1}{24} c_{22} = 0.
\end{equation}
\begin{equation}\label{eq21}
 \< \d^j \d^\nu \d^\nu \> \< \d^x \d^i \d^\alpha \d^\beta \> \< \d^\alpha \>_1 \< \d^\beta \>_1 
: \quad
\frac{13}{10} c_3 + c_{15} - c_{19} + \frac{1}{8} c_{29}=0.
\end{equation}
\begin{equation}\label{eq22}
 \< \d^j \d^x \d^\mu \d^\mu \> \< \d^i \d^\alpha \d^\beta \> \< \d^\alpha \>_1 \< \d^\beta \>_1 
: \quad
\frac{7}{10} c_4 + c_{17} - c_{19} + \frac{1}{24} c_{28} = 0.
\end{equation}
\begin{equation}\label{eq23}
 \< \d^j \d^\nu \d^\nu \> \< \d^x \d^i \d^\alpha \d^\alpha \d^\beta \> \< \d^\beta \>_1 
: \quad
\frac{23}{240} c_3 + c_8 - 2 c_{10} + \frac{1}{24} c_{15}
+ \frac{1}{12} c_{19} = 0.
\end{equation}
\begin{equation}\label{eq24}
 \< \d^j \d^x \d^\mu \d^\mu \> \< \d^i \d^\alpha \d^\alpha \d^\beta \> \< \d^\beta \>_1 
: \quad
\frac{13}{240} c_4 + c_7 - 2 c_{10} + \frac{1}{24} c_{17}
+ \frac{1}{24} c_{18} = 0.
\end{equation}
\begin{equation}\label{eq25}
 \< \d^i \d^\alpha \d^\alpha \d^\beta \d^\beta \> \< \d^x \d^j \d^\mu \> \< \d^\mu \>_1 
: \quad
\frac{1}{960} c_6 + c_9 - c_{10} + \frac{1}{24} c_{16} = 0.
\end{equation}
\begin{equation}\label{eq35}
\< \d^x \d^\mu \d^j \> \< \d^i \d^\alpha \d^\alpha \d^\nu \> \< \d^\mu \>_1 \< \d^\nu \>_1 
: \quad
\frac{13}{240} c_6 + c_{16} + c_{18} - 2 c_{19} + \frac{1}{24} c_{27} = 0.
\end{equation}
\begin{equation}\label{eq36}
 \< \d^j \d^\nu \d^\nu \> \< \d^x \d^i \d^\alpha \d^\alpha \d^\mu \d^\mu \> 
: \quad
\frac{1}{576} c_3 + \frac{1}{24} c_8 + \frac{1}{24} c_{10} - 3 c_{30} = 0.
\end{equation}
\begin{equation}\label{eq37}
 \< \d^j \d^x \d^\mu \d^\mu \> \< \d^i \d^\alpha \d^\alpha \d^\nu \d^\nu \> 
: \quad
\frac{1}{960} c_4 + \frac{1}{24} c_7 + \frac{1}{24} c_9 - 3 c_{30} = 0.
\end{equation}

There are several terms of the form $* \< \d^i \>_1 $.
If we remove the $ \< \d^i \>_1 $, they become terms in $\ocM_{2,2}$ of
codimension $3$, and there is a relation between them which we can find by
using Getzler's relation (with a two descendences on $x$ and one descendence
on $j$).
The relation is
\begin{eqnarray*}
0 & = &
- \frac{3}{40} \< \d^x \d^\mu \d^\alpha \d^\alpha \> \< \d^y \d^\mu \d^\nu \> \< \d^\nu \>_1
+ \frac{3}{40} \< \d^x \d^\mu \d^\nu \d^\alpha \> \< \d^\alpha \d^y \d^\mu \> \< \d^\nu \>_1 \\
& & - \frac{7}{120} \< \d^x \d^\mu \d^\alpha \> \< \d^\alpha \d^y \d^\mu \d^\nu \> \< \d^\nu \>_1
+ \frac{7}{120} \< \d^x \d^\nu \d^\alpha \> \< \d^\alpha \d^y \d^\mu \d^\mu \> \< \d^\nu \>_1 \\
& & + \frac{1}{120} \< \d^x \d^\mu \d^\alpha \> \< \d^\alpha \d^\mu \d^\nu \> \< \d^y \d^\nu \>_1
- \frac{1}{120} \< \d^x \d^y \d^\nu \d^\alpha \> \< \d^\alpha \>_1 \< \d^\mu \d^\mu \d^\nu \> \\
& & - \frac{1}{120} \< \d^x \d^y \d^\alpha \> \< \d^\nu \d^\alpha \>_1 \< \d^\mu \d^\mu \d^\nu \>
+ \frac{1}{120} \< \d^y \d^\mu \d^\alpha \> \< \d^x \d^\alpha \>_1 \< \d^\mu \d^\nu \d^\nu \> \\
& & + \textrm{other terms with all vertices of genus $0$} .
\end{eqnarray*}
We are going to solve this relation for the term
$ \< \d^x \d^\mu \d^\alpha \> \< \d^\alpha \d^\mu \d^\nu \> \< \d^j \d^\nu \>_1 \< \d^i \>_1 $
and find an equation for all of the other terms.
Among those, seven of them are of the form $ * \< \d^i \>_1 \< \d^\mu \>_1 $
and are related by WDVV. They can be written in terms of the following $4$ 
independent vectors.
\begin{equation}\label{eq30}
\< \d^\alpha \d^\alpha \d^j \d^\nu \> \< \d^\nu \d^x \d^\mu \> \< \d^i \>_1 \< \d^\mu \>_1 
: \quad 
- \frac{3}{20} c_1 - c_2 - 2 c_3 - 2 c_5 - \frac{1}{24} c_6 - c_{16} + 2 c_{19} - c_{24} + \frac{1}{24} c_{27} = 0.
\end{equation}
\begin{equation}\label{eq31}
\< \d^\alpha \d^\alpha \d^x \d^\nu \> \< \d^\nu \d^j \d^\mu \> \< \d^i \>_1 \< \d^\mu \>_1 
: \quad 
- c_2 - c_4 = 0.
\end{equation}
\begin{equation}\label{eq32}
\< \d^\alpha \d^\alpha \d^\mu \d^\nu \> \< \d^\nu \d^x \d^j \> \< \d^i \>_1 \< \d^\mu \>_1 
: \quad 
\frac{11}{240} c_1 + c_2 + 2 c_3 + c_5 - c_{18} + c_{24} = 0.
\end{equation}
\begin{equation}\label{eq33}
\< \d^\alpha \d^j \d^\nu \> \< \d^\nu \d^x \d^\mu \d^\alpha \> \< \d^i \>_1 \< \d^\mu \>_1 
: \quad 
- \frac{1}{10} c_1 - 2 c_3 - c_5 + 4 c_{19} - c_{23} - c_{24} = 0.
\end{equation}
There are additional $4$ independent vectors
\begin{equation}\label{eq26}
 \< \d^x \d^\beta \d^\gamma \> \< \d^\gamma \d^j \d^\alpha \> \< \d^\alpha \d^\beta \>_1 \< \d^i \>_1 
: \quad
\frac{1}{10} c_1 + c_5 + c_{22} - c_{23} - 2 c_{25} + 2 c_{26} = 0.
\end{equation}
\begin{equation}\label{eq27}
\< \d^x \d^\beta \d^\mu \> \< \d^\mu \d^j \d^\alpha \> \< \d^\alpha \>_1 \< \d^\beta \>_1 \< \d^i \>_1 
: \quad
- \frac{7}{10} c_1 - c_6 - c_{28} + 3 c_{29} = 0.
\end{equation}
\begin{equation}\label{eq28}
 \< \d^\mu \d^\nu \d^\nu \> \< \d^i \d^\mu \d^\alpha \> \< \d^x \d^\alpha \>_1 \< \d^j \>_1 
: \quad
\frac{1}{20} c_1 + c_2 + c_3 + c_5 - c_{14} + c_{15} + c_{18} - c_{22} = 0.
\end{equation}
\begin{equation}\label{eq29}
 \< \d^\alpha \d^\beta \d^\beta \> \< \d^x \d^j \d^\gamma \> \< \d^\alpha \d^\gamma \>_1 \< \d^i \>_1 
: \quad
- \frac{11}{240} c_1 - c_2 - 2 c_3 - c_5 - c_{14} + c_{15} = 0.
\end{equation}

The remaining connected strata are all of codimension $3$ in $\ocM_{2,3}$.
There are four induced equations.
First of all, by taking Getzler's relation in $\ocM_{1,4}$, adding another
marked point, and then identifying either two of the first four marked points
or the fifth marked point with one of the others.
Secondly, by taking the Belorousski-Pandharipande relation, adding a
descendence at the marked point $x$, and then simplifying.
The relations we obtain are the following:
\begin{eqnarray*}
0 & = & \< \d^x \d^i \d^\mu \> \< \d^j \d^\alpha \d^\nu \> \< \d^\alpha \d^\mu \d^\nu \>_1
+ \< \d^x \d^j \d^\mu \> \< \d^i \d^\alpha \d^\nu \> \< \d^\alpha \d^\mu \d^\nu \>_1 \\ & & \hspace{.5cm}
+ \< \d^x \d^\alpha \d^\mu \> \< \d^i \d^j \d^\nu \> \< \d^\alpha \d^\mu \d^\nu \>_1
- \< \d^x \d^i \d^\mu \> \< \d^j \d^\mu \d^\nu \> \< \d^\alpha \d^\alpha \d^\nu \>_1 \\ & & \hspace{.5cm}
- \< \d^x \d^i \d^\mu \> \< \d^\alpha \d^\mu \d^\nu \> \< \d^j \d^\alpha \d^\nu \>_1
- \< \d^x \d^j \d^\mu \> \< \d^\alpha \d^\mu \d^\nu \> \< \d^i \d^\alpha \d^\nu \>_1 \\ & & \hspace{.5cm}
- \< \d^i \d^j \d^\mu \> \< \d^\alpha \d^\mu \d^\nu \> \< \d^x \d^\alpha \d^\nu \>_1
+ \textrm{other terms},
\end{eqnarray*}
\begin{eqnarray*}
0 & = & \< \d^x \d^i \d^\mu \> \< \d^\alpha \d^\alpha \d^\nu \> \< \d^j \d^\mu \d^\nu \>_1
- \< \d^x \d^\alpha \d^\mu \> \< \d^\alpha \d^\mu \d^\nu \> \< \d^i \d^j \d^\nu \>_1 \\ & & \hspace{.5cm}
- \< \d^i \d^\alpha \d^\mu \> \< \d^\alpha \d^\mu \d^\nu \> \< \d^x \d^j \d^\nu \>_1
+ \textrm{other terms},
\end{eqnarray*}
\begin{eqnarray*}
0 & = & \< \d^i \d^j \d^\mu \> \< \d^\alpha \d^\alpha \d^\nu \> \< \d^x \d^\mu \d^\nu \>_1
- \< \d^i \d^\alpha \d^\mu \> \< \d^\alpha \d^\mu \d^\nu \> \< \d^j \d^x \d^\nu \>_1 \\ & & \hspace{.5cm}
- \< \d^j \d^\alpha \d^\mu \> \< \d^\alpha \d^\mu \d^\nu \> \< \d^i \d^x \d^\nu \>_1
+ \textrm{other terms},
\end{eqnarray*}
\begin{eqnarray*}
0 & = & - \< \d^\mu \d^\alpha \d^\nu \> \< \d^\mu \d^i \d^j \> \< \d^\nu \>_1 \< \d^x \d^\alpha \>_1
- \< \d^x \d^j \d^\alpha \> \< \d^\alpha \d^\mu \d^\nu \> \< \d^\mu \>_1 \< \d^\nu \d^i \>_1 \\
& & - \< \d^x \d^i \d^\alpha \> \< \d^\alpha \d^\mu \d^\nu \> \< \d^\mu \>_1 \< \d^\nu \d^j \>_1
- \< \d^x \d^\mu \d^\alpha \> \< \d^\mu \d^i \d^j \> \< \d^\alpha \d^\nu \>_1 \< \d^\nu \>_1 \\
& & + \< \d^x \d^\alpha \d^\nu \> \< \d^\mu \d^i \d^j \> \< \d^\alpha \d^\mu \>_1 \< \d^\nu \>_1
+ \< \d^x \d^j \d^\alpha \> \< \d^i \d^\mu \d^\nu \> \< \d^\mu \>_1 \< \d^\nu \d^\alpha \>_1 \\
& & + \< \d^x \d^i \d^\alpha \> \< \d^j \d^\mu \d^\nu \> \< \d^\mu \>_1 \< \d^\nu \d^\alpha \>_1
+ \textrm{other terms}.
\end{eqnarray*}
After solving the four relations above for the terms
$ \< \d^x \d^i \d^\mu \> \< \d^j \d^\mu \d^\nu \> \< \d^\alpha \d^\alpha \d^\nu \>_1, $ 
$ \< \d^x \d^i \d^\mu \> \< \d^\alpha \d^\alpha \d^\nu \> \< \d^j \d^\mu \d^\nu \>_1, $
$ \< \d^i \d^j \d^\mu \> \< \d^\alpha \d^\alpha \d^\nu \> \< \d^x \d^\mu \d^\nu \>_1, $
and
$ \< \d^x \d^\mu \d^\alpha \> \< \d^\mu \d^i \d^j \> \< \d^\alpha \d^\nu \>_1 \< \d^\nu \>_1, $
we obtain the following equations for the other terms.
\begin{equation}\label{eq38}
\< \d^x \d^\alpha \d^\beta \> \< \d^\beta \d^i \d^j \> \< \d^\alpha_1 \>_2
: \quad
- \frac{1}{2} c_1 + 4 c_5 - \frac{1}{2} c_6 = 0
\end{equation}
\begin{equation}\label{eq39}
\< \d^\alpha \d^\beta \d^\beta \> \< \d^x \d^\alpha \d^\gamma \>\< \d^i \d^j \d^\gamma \>_1
: \quad
- \frac{1}{40} c_1 - \frac{1}{2} c_2 - \frac{1}{2} c_3 - \frac{1}{2} c_5 + 2 c_8 - \frac{1}{2} c_{15} = 0.
\end{equation}
\begin{equation}\label{eq40}
\< \d^\mu \d^\nu \d^\nu \> \< \d^i \d^\mu \d^\alpha \> \< \d^j \d^x \d^\alpha \>_1
: \quad
\frac{1}{240} c_1 - c_3 - \frac{1}{60} c_5 + 4 c_8 + 2 c_{12} - c_{15} - 3 c_{21} = 0.
\end{equation}
\begin{equation}\label{eq44}
 \< \d^x \d^\beta \d^\gamma \> \< \d^\gamma \d^i \d^\alpha \> \< \d^j \d^\alpha \d^\beta \>_1 
: \quad
- \frac{1}{10} c_1 - \frac{29}{30} c_5 - 4 c_{11} + c_{16} + 3 c_{20} - 3 c_{21} = 0.
\end{equation}
\begin{equation}\label{eq43}
 \< \d^i \d^j \d^\beta \> \< \d^\beta \d^\mu \d^\nu \> \< \d^x \d^\mu \d^\nu \>_1 
: \quad
\frac{1}{30} c_5 - 2 c_{11} + \frac{1}{2} c_{16} + 6 c_{21} - \frac{1}{2} c_{24} = 0.
\end{equation}
\begin{equation}\label{eq41}
 \< \d^x \d^j \d^\nu \> \< \d^i \d^\alpha \d^\beta \> \< \d^\nu \d^\alpha \d^\beta \>_1 
: \quad
\frac{1}{15} c_5 + 8 c_{11} - c_{16} - 3 c_{20} + 3 c_{21} = 0.
\end{equation}
\begin{equation}\label{eq42}
 \< \d^x \d^\mu \d^\nu \> \< \d^i \d^j \d^\beta \> \< \d^\mu \d^\nu \d^\beta \>_1 
: \quad
- \frac{1}{30} c_5 + 4 c_{11} - \frac{1}{2} c_{16} - \frac{1}{2} c_{25} = 0.
\end{equation}
\begin{equation}\label{eq46}
 \< \d^x \d^\beta \d^\gamma \> \< \d^\gamma \d^i \d^\alpha \> \< \d^j \d^\alpha \>_1 \< \d^\beta \>_1 
: \quad
- \frac{7}{5} c_1 + \frac{7}{5} c_5 - c_6 + 2 c_{23} - 2 c_{24} - 6 c_{25} + 2 c_{26} + c_{27} = 0.
\end{equation}
\begin{equation}\label{eq45}
 \< \d^i \d^j \d^\beta \> \< \d^\beta \d^\mu \d^\nu \> \< \d^x \d^\mu \>_1 
: \quad
\frac{4}{5} c_5 + 2 c_{22} + 2 c_{24} - 2 c_{25} + \frac{1}{2} c_{27} - c_{28} = 0.
\end{equation}
\begin{equation}\label{eq48}
 \< \d^x \d^\mu \d^j \> \< \d^i \d^\alpha \d^\beta \> \< \d^\mu \d^\alpha \>_1 \< \d^\beta \>_1 
: \quad
\frac{3}{5} c_5 - 2 c_{23} + 2 c_{24} + 6 c_{25} + 2 c_{26} - c_{27} = 0.
\end{equation}
\begin{equation}\label{eq47}
 \< \d^x \d^\mu \d^\nu \> \< \d^i \d^j \d^\beta \> \< \d^\nu \d^\beta \>_1 \< \d^\mu \>_1 
: \quad
- \frac{4}{5} c_5 + 2 c_{16} + 2 c_{25} - c_{27} = 0.
\end{equation}

Solving the equations (\ref{eq1})-(\ref{eq47}) gives the following
coefficients (let $c_1= - 1$):
$$ \begin{matrix}
& c_2 = - \displaystyle \frac{5}{72} 
& c_3 = - \displaystyle \frac{1}{252}
& c_4 = \displaystyle \frac{5}{72}
& c_5 = \displaystyle \frac{5}{42}
& c_6 = \displaystyle \frac{41}{21} \\ \\
& c_7 = - \displaystyle \frac{11}{40320}
& c_8 = - \displaystyle \frac{1}{13440} 
& c_9 = - \displaystyle \frac{1}{8064}
& c_{10} = \displaystyle \frac{191}{120960}
& c_{11} = - \displaystyle \frac{1}{5040} \\ \\ 
& c_{12} = - \displaystyle \frac{1}{4032} 
& c_{13} = \displaystyle \frac{17}{2880}
& c_{14} = \displaystyle \frac{1}{840}
& c_{15} = - \displaystyle \frac{1}{336}
& c_{16} = - \displaystyle \frac{1}{126} \\ \\
& c_{17} = - \displaystyle \frac{23}{5040}
& c_{18} = - \displaystyle \frac{17}{5040}
& c_{19} = \displaystyle \frac{113}{2520}
& c_{20} = \displaystyle \frac{1}{210}
& c_{21} = 0 \\ \\
& c_{22} = - \displaystyle \frac{1}{84}
& c_{23} = \displaystyle \frac{211}{1260}
& c_{24} = \displaystyle \frac{1}{1260}
& c_{25} = - \displaystyle \frac{1}{630}
& c_{26} = \displaystyle \frac{11}{140} \\ \\
& c_{27} = - \displaystyle \frac{4}{35}
& c_{28} = \displaystyle \frac{2}{105}
& c_{29} = \displaystyle \frac{89}{210}
& c_{30} = \displaystyle \frac{1}{53760}
\end{matrix} $$


\end{document}